\documentclass[graybox]{svmult}

\usepackage{type1cm}        % activate if the above 3 fonts are
\usepackage{makeidx}         % allows index generation
\usepackage{graphicx}        % standard LaTeX graphics tool
\usepackage{multicol}        % used for the two-column index
\usepackage[bottom]{footmisc}% places footnotes at page bottom

\usepackage{newtxtext}       % 
\usepackage[varvw]{newtxmath}       % selects Times Roman as basic font
\usepackage{bm}
\usepackage{empheq}

\usepackage{subfig}
\usepackage[percent]{overpic}
\usepackage{tikz}
\usepackage{circuitikz}
\usetikzlibrary{calc,fit,backgrounds}
\usetikzlibrary{shapes.geometric, arrows}
\usetikzlibrary{chains,
                fit,
                positioning,
                shapes}
\usetikzlibrary{decorations.pathreplacing}
\usepackage{color}
\usepackage{pict2e}
\usepackage[numbers]{natbib}
\usepackage{xspace}
\usepackage[dvipsnames,table]{xcolor}
\makeindex             % used for the subject index
\begin{document}

%\title*{Reduced Order Models in Cardiovascular Flow Analysis: Methods and Applications}

\title*{A Review of Equation-Based and Data-Driven Reduced Order Models featuring a Hybrid cardiovascular application}
\titlerunning{A review of ROMs featuring a hybrid cardiovascular application}
\authorrunning{P. Siena, P. C. Africa, M. Girfoglio, G. Rozza}

\author{Pierfrancesco Siena\orcidID{0009-0007-2124-0598} and\\ Pasquale Claudio Africa \orcidID{0000-0002-0706-8564} and\\ Michele Girfoglio \orcidID{0000-0003-1766-2265} and\\ Gianluigi Rozza \orcidID{0000-0002-0810-8812}}
\institute{Pierfrancesco Siena \at SISSA, International School for Advanced Studies, Mathematics Area, mathLab, via Bonomea 265, I-34136 Trieste, Italy \email{psiena@sissa.it}
\and Pasquale Claudio Africa \at SISSA, International School for Advanced Studies, Mathematics Area, mathLab, via Bonomea 265, I-34136 Trieste, Italy \email{pafrica@sissa.it}\and Michele Girfoglio \at University of Palermo, Department of Engineering,  Viale delle Scienze, Ed. 7, 90128 Palermo, Italy  \email{michele.girfoglio@unipa.it}\and Gianluigi Rozza \at SISSA, International School for Advanced Studies, Mathematics Area, mathLab, via Bonomea 265, I-34136 Trieste, Italy \email{grozza@sissa.it}}
\maketitle

\abstract*{Each chapter should be preceded by an abstract (no more than 200 words) that summarizes the content. The abstract will appear \textit{online} at \url{www.SpringerLink.com} and be available with unrestricted access. This allows unregistered users to read the abstract as a teaser for the complete chapter.
Please use the 'starred' version of the \texttt{abstract} command for typesetting the text of the online abstracts (cf. source file of this chapter template \texttt{abstract}) and include them with the source files of your manuscript. Use the plain \texttt{abstract} command if the abstract is also to appear in the printed version of the book.}
\abstract{Cardiovascular diseases are a leading cause of death in the world, driving the development of patient-specific and benchmark models for blood flow analysis.
This chapter provides a theoretical overview of the main categories of Reduced Order Models (ROMs), focusing on both projection-based and data-driven approaches within a classical setup. We then present a hybrid ROM tailored for simulating blood flow in a patient-specific aortic geometry. The proposed methodology integrates projection-based techniques with neural network-enhanced data-driven components, incorporating a lifting function strategy to enforce physiologically realistic outflow pressure conditions. This hybrid methodology enables a substantial reduction in computational cost while mantaining high fidelity in reconstructing both velocity and pressure fields. We compare the full- and reduced-order solutions in details and critically assess the advantages and limitations of ROMs in patient-specific cardiovascular modeling.}

\section*{Introduction}
\label{sec:introduction}
\textcolor{black}{Cardiovascular diseases are among the leading causes of death worldwide, posing a significant burden on global healthcare systems. Conditions such as coronary artery disease, stroke and heart failure demand more effective diagnostic and therapeutic strategies. In recent years, computational modeling has gained prominence as a powerful approach to support cardiovascular research and clinical decision-making. By integrating mathematical formulations with physiological and clinical data, these models allow for detailed simulations of cardiovascular dynamics, including blood flow, offering valuable insights into disease mechanisms and treatment outcomes.}

In computational modeling, achieving a balance between efficiency and accuracy remains a central challenge for researchers and engineers, \textcolor{black}{particularly when addressing real-world problems.} While high-fidelity discretization methods, such as finite element and finite volume, enable detailed simulations of complex physical phenomena, their computational demands can be prohibitive, \textcolor{black}{due to the challenging nature of the domain or of the flow dynamics}. These Full Order Models (FOMs) require significant time and computational resources, making them impractical for real-world applications, particularly when many physical and/or geometric parameters are involved.

To address these limitations, Reduced Order Models (ROMs) provide an efficient alternative, significantly reducing simulation time. While this computational speedup is highly beneficial, it often raises concerns about the accuracy of the resulting solutions. The capability of ROMs to achieve real-time evaluations is closely related to the fact that the parameters involved in the problem usually operate in a much lower-dimensional space than the number of degrees of freedom associated with the FOM. ROMs operate within a completely decoupled \emph{offline}-\emph{online} paradigm \cite{benner2020model,rozza20201}, where the computationally expensive tasks are performed in an offline phase, enabling rapid and efficient evaluations during the online stage. More in detail, during the \emph{offline} stage, the FOM is solved for specific parameter values; afterward, a compact reduced basis is constructed from the database of FOM solutions. In the \emph{online} stage, the reduced coefficients are typically predicted using interpolation or regression techniques (for data-driven ROMs), or by solving a system of algebraic or ordinary differential equations (for projection-based ROMs).

This chapter begins by examining the classical framework of ROMs through a simple, well-known problem that is often used as a foundational example to illustrate the key features of a general ROM setup. Although basic, this formulation shares many fundamental principles with more complex problems, making it a valuable reference.
We first provide a detailed introduction to projection-based (or equation-based) techniques, which reduce the dimensionality of a FOM by projecting its governing equations onto a lower-dimensional subspace. Next, we explore data-driven approaches, which rely on mathematical and machine learning techniques to approximate system dynamics directly from data, without explicit knowledge of the governing equations. All these techniques are provided in Sec. \ref{sec:rom}.

In Sec. \ref{sec:aorta}, previous foundations are combined to address a cardiovascular application through a hybrid approach, which combine the robustness of equation-based reduction with the flexibility and adaptability of data-driven models. %Hybrid methods are useful for example when the underlying physical model is only partially known or too costly to evaluate. %In such cases, the main dynamics are typically described by reduced equations, while specific quantities of interest, such as nonlinear source terms, stabilization terms, effective parameters, or closure models, can be learned from data. These components can be generalized and approximated using interpolation, regression or machine learning techniques.

%\textcolor{black}{This chapter examines the features of ROMs, both the accuracy and efficiency, emphasizing their practical impact in the numerical simulations of the cardiovascular system. %The goal is not only to assess their computational speed but also to evaluate their accuracy and robustness in real-world applications.
%We consider ROMs based on three different strategies: (i) data-driven methods, which leverage mathematical and machine learning techniques to approximate system dynamics directly from data; (ii) projection-based techniques, which reduce the dimensionality of the problem by projecting the FOM equations onto a lower-dimensional subspace; and (iii) hybrid approaches, which combine elements of both data-driven and projection-based methods.}% to enhance accuracy and efficiency. 

%The chapter is structured as follows. Sec. \ref{sec:rom} is dedicated to the mathematical description of techniques used for projection-based and data-driven ROMs. Then, a patient-specific study is shown in Sec. \ref{sec:aorta} for the study of the hemodynamics of the aortic arch. %and an idealized biomedical device, the Food and Drug Administration benchmark (FDA), is proposed in Sec. \ref{sec:fda}.

\section{Reduced Order Models}
\label{sec:rom}
In this section, we provide a mathematical overview of ROMs. We focus on the two main categories of ROMs \cite{aromabook,benner2020model,benner2020model2,degruyter1,rozza20201}: projection-based ROMs, which directly leverage the governing equations of the problem (see, e.g., \cite{stabile2017pod,stabile2018finite,girfoglio2021pod,girfoglio2021pressure,pichi2022driving,pichi2019reduced,pichi2020reduced}), and data-driven ROMs, which rely solely on data without explicitly incorporating the governing equations (see, e.g., \cite{siena2023fast,Siena2022,balzotti2023reduced,zainib2022chapter,balzottidata2022,papapicco2021neural,hijazi2020data,shah2021finite,coscia2023physics,pichi2023artificial}).

In the case of projection-based ROMs, we focus on Galerkin projection methods, where the full-order governing equations are projected onto a reduced basis space, yielding a system of differential-algebraic equations whose unknowns are the reduced coefficients \cite{rowley2004model,iollo2000stability}.

To formalize the problem, we now introduce key concepts essential for the development of the methods presented in the following sections. 

Let us define the computational domain $\Omega \subset \mathbf{R}^d$ where $d=1,2,3$, with boundary $\partial \Omega$. We denote by $\Gamma^D_i$ the regions of $\partial \Omega$ where Dirichlet boundary conditions are enforced for $1\leq i \leq d_v$. Let $w:\Omega \mapsto \mathbf{R}^{d_v}$ represent a vector field ($d_v = d$) or a scalar field ($d_v = 1$).

The following scalar space is defined along with an appropriate norm $\lVert \cdot \rVert_{\mathbf{V}}$:
\begin{equation}
    \mathbf{V}_i(\Omega)=\{ v\in H^1(\Omega) : v_{\mid \Gamma^D_i}=0 \}, \qquad  1\leq i \leq d_v,
\end{equation}
where $H^1_0(\Omega) \subset \mathbf{V}_i(\Omega) \subset H^1(\Omega)$. To proceed with the discretization, we consider a finite-dimensional subspace $\mathbf{V}_{\delta} \subset \mathbf{V}$ with size $\text{dim}(\mathbf{V}_{\delta}) = N_\delta$. 

To support the analysis of parametric cases, we define a closed parameter domain $\mathbf{P} \subset \mathbf{R}^p$. A general element of $\mathbf{P}$ is denoted by $\mu = (\mu_1, \dots, \mu_p)$, and the following parametric field variable for a parameter value is introduced  $u(\mu)=(u_1 (\mu),\dots, u_{d_v}(\mu)):\mathbf{P}\mapsto\mathbf{V}$. The discretized form of the parameter space is identified by $\mathbf{P}_h \subset \mathbf{P}$. 

Finally, let us introduce the solution manifold:
\begin{equation}
    \mathcal{M} = \{u(\mu) : \mu \in \mathbf{P} \} \subset \mathbf{V},
\end{equation}
collecting all the exact solutions $u(\mu)\in\mathbf{V}$ varying the parameter $\mu$.
Similarly, at a discrete level, the solution manifold is:
\begin{equation}
    \mathcal{M}_{\delta}=\{ u_{\delta}(\mu) : \mu \in \mathbf{P}\} \subset \mathbf{V}_\delta,
\end{equation}
where $u_{\delta}(\mu) \in \mathbf{V}_{\delta}$ are the high fidelity solutions, computed with finite element or finite volume methods.

\subsection{Building reduced basis space}

To build a reduced basis space, the Proper Orthogonal Decomposition (POD) is among the most commonly employed techniques. Various other approaches, both linear and nonlinear, have been proposed in the literature, including the greedy algorithm, Proper Generalized Decomposition (PGD), factor analysis, independent component analysis, and autoencoders \cite{cunningham2015linear,fu2021data,lee2020model,gonzalez2018deep,kashima2016nonlinear}. However, in this work, we focus exclusively on POD, as it is the method utilized in the application discussed in the following.

\subsubsection{Proper orthogonal decomposition}
\label{sec:pod}
%POD is a method suitable for data compression \cite{eckart1936approximation}. Once the parameter space is discretized and FOM solutions are collected for each element of $\mathbf{P}_h$, POD allows to retain only the essential information about the system at hand. %then a compact formulation of the problem can be achieved. The $N$-dimensional POD space can be computed as the solution to the following minimization problem:

The POD is an effective technique for data compression \cite{eckart1936approximation}. %After discretizing the parameter space and collecting the FOM solutions for each element of $\mathbf{P}_h$, POD 
It enables the extraction of the most relevant features of the system while discarding redundant information. The POD algorithm leads to a compact and accurate representation of the problem. The resulting $N$-dimensional space is obtained by solving the following minimization problem:
\begin{equation}
    \min_{\mathbf{V}_{\text{rb}}:|\mathbf{V}_{\text{rb}}|=N} \Bigg( \int_{\mu \in \mathbf{P}} \inf_{v_{\text{rb}}\in \mathbf{V}_{\text{rb}}} \lVert u_{\delta}(\mu)-v_{\text{rb}} \rVert_{\mathbf{V}}^2 d\mu \Bigg)^{1/2}. \label{eq:POD1}
\end{equation}
In discrete form eq. \eqref{eq:POD1} yields:
\begin{equation}
     \min_{\mathbf{V}_{\text{rb}}:|\mathbf{V}_{\text{rb}}|=N} \Bigg( \frac{1}{M}\sum_{\mu \in \mathbf{P}_h} \inf_{v_{\text{rb}}\in \mathbf{V}_{\text{rb}}} \lVert u_{\delta}(\mu)-v_{\text{rb}} \rVert_{\mathbf{V}}^2  \Bigg)^{1/2}.
     \label{opt-crit}
\end{equation}
%Let sort the elements of $\mathbf{P}_h$
Let the elements of $\mathbf{P}_h$ be 
\begin{equation}
    \{\mu_1,\dots,\mu_M\},
\end{equation}
and denote by
\begin{equation}
    \{ \psi_1, \dots, \psi_M \},
\end{equation}
the elements of $\mathcal{M}_{\delta}(\mathbf{P}_h)$, $\psi_m=u_{\delta}(\mu_m)$ for $m=1,\dots,M$ (\textit{i.e.}, the high-fidelity snapshots). Furthermore, let us define $\mathbf{V}_{\mathcal{M}}=\text{span}\{ u_{\delta}(\mu) : \mu \in \mathbf{P}_h \}$ and the symmetric and linear operator $\mathcal{C}:\mathbf{V}_{\mathcal{M}}\mapsto\mathbf{V}_{\mathcal{M}}$ such that:
\begin{equation}
    \mathcal{C}(v_{\delta})=\frac{1}{M}\sum_{m=1}^M (v_{\delta},\psi_m)_{\mathbf{V}}\,\psi_m, \quad v_{\delta} \in \mathbf{V}_{\mathcal{M}}.
    \label{linear-op}
\end{equation}
Its eigenvalues and normalized eigenvectors are identified as $(\lambda_i,\xi_i)\in \mathbf{R}\times \mathbf{\mathbf{V}_{\mathcal{M}}}$ and fulfill:
\begin{equation}
    (\mathcal{C}(\xi_i),\psi_m)_{\mathbf{V}}=\lambda_i(\xi_i,\psi_m)_{\mathbf{V}}, \quad 1 \le m \le M,
    \label{eig-p}
\end{equation}
with 
$$\lambda_1 \ge \lambda_2 \ge \dots \ge \lambda_M.$$
%Here we assume that the eigenvalues are sorted in descending order, $\lambda_1 \ge \lambda_2 \ge \dots \ge \lambda_M$. 
The eigenvectors $\{\xi_1,\dots,\xi_M\}$ represent the basis functions for the space $\mathbf{V}_{\mathcal{M}}$ and the first $N \ll M$ eigenfunctions $\{\xi_1,\dots,\xi_N\}$ yield the $N$-dimensional reduced space
\begin{equation}
    \mathbf{V}_{\text{POD}} = \text{span} \{\xi_1,\dots,\xi_N \}.
\end{equation}
%$\mathbf{V}_{\text{POD}} = \text{span} \{\xi_1,\dots,\xi_N \}$. %which fulfills {\color{cyan} Eq.} \eqref{opt-crit}.

The error introduced by replacing the elements of $\mathcal{M}_{\delta}(\mathbf{P}_h)$ by $\mathbf{V}_{\text{POD}}$  is given by the sum of the neglected eigenvalues \cite{quarteroni2015reduced}:
%By projecting the elements of $\mathcal{M}_{\delta}(\mathbf{P}_h)$ onto $\mathbf{V}_{\text{POD}}$, it can be proven that the error is related to the neglected eigenvalues:
\begin{equation}
    \frac{1}{M}\sum_{m=1}^M\lVert \psi_m - P_N[\psi_m] \rVert_{\mathbf{V}}^2 = \sum_{m=N+1}^M \lambda_m,
\end{equation}
where $P_N[\psi_m] = \sum_{i=1}^N(\psi_m,\xi_i)_{\mathbf{V}}\,\xi_i$ is the projection of $\psi_m$ onto $\mathbf{V}_{\text{POD}}$.
As a consequence of the orthonormality of the eigenvectors in $\lVert \cdot \rVert_{\ell^2(\mathbf{R}^M)}$, we have:
\begin{equation}
    (\xi_m,\xi_q)_{\mathbf{V}}=M\lambda_i\delta_{mq}, \quad 1\le m,q \le M,
\end{equation}
where $\delta_{mq}$ is the Kronecker delta.

The correlation matrix $C\in\mathbf{R}^{M\times M}$ representing the linear operator \eqref{linear-op} can be expressed as:
\begin{equation}
    C_{m,q}=\frac{1}{M}(\psi_m,\psi_q)_{\mathbf{V}}, \quad 1 \le m,q \le M.
\end{equation}
Therefore the eigenvalues problem \eqref{eig-p} is reformulated as:
\begin{equation}
    C v_i = \lambda_i v_i, \quad 1 \le i \le N.
\end{equation}
Finally, the orthogonal basis functions are given by:
\begin{equation}
    \xi_i = \frac{1}{\sqrt{M}}\sum_{m=1}^M (v_i)_m \psi_m, \quad 1 \le i \le N,
\end{equation}
where $(v_i)_m$ denotes the $m$-th element of the eigenvector $v_i \in \mathbf{R}^M$ .

The computational cost to implement the POD algorithm can be very large, because it is impossible to know how many high-fidelity solutions are needed to guarantee a good approximation of the system behavior, and it depends on the problem at hand. Typically $M\gg N$ instances of the FOM need to be solved during the offline stage and, when $M$ and $N_{\delta}$ are large, the cost of extracting the eigenfunctions rises, as it scales as $\mathcal{O}(NN_{\delta}^2)$.

\subsection{Projection-based models}
\label{sec:proj-based}

To simplify the discussion \textcolor{black}{while preserving core principles}, we focus on a stationary problem and adopt a finite element formulation for its discretization. Given $\mu \in \mathbf{P} $, the abstract formulation of a stationary problem evaluates $s(\mu) = l(u(\mu);\mu)$ where $u(\mu)\in \mathbf{V}$ satisfies:
\begin{equation}
      a(u(\mu),v;\mu) = f(v;\mu) \qquad \forall v \in \mathbf{V}.
      \label{varitional-problem}
\end{equation}
The form $a:\mathbf{V}\times\mathbf{V}\times\mathbf{P}\mapsto\mathbf{R}$ is bilinear with respect to $\mathbf{V} \times \mathbf{V}$, $f:\mathbf{V}\times\mathbf{P}\mapsto\mathbf{R}$ and $l:\mathbf{V}\times\mathbf{P}\mapsto\mathbf{R}$ are linear with respect to $\mathbf{V}$, and $s:\mathbf{P}\mapsto\mathbf{R}$ is the scalar output of the model. The problem \eqref{varitional-problem} represents our FOM in its continuous formulation. Throughout this analysis, we assume that the compliance hypothesis holds. Specifically, we consider:
\begin{itemize}
    \item $l(\cdot;\mu) = f(\cdot;\mu) \quad \forall \mu \in \mathbf{P}$;
    \item $a(\cdot,\cdot;\mu)$ symmetric $\forall \mu \in \mathbf{P}$.
\end{itemize}
This assumption significantly streamlines the presentation while still capturing the essential ROM principles. Furthermore, we assume the well-posedness of problem \eqref{varitional-problem} and we observe that the Lax-Milgram theorem \cite{quarteroni2008numerical} holds. %states that if 
% \begin{itemize}
%     \item the bilinear form $a(\cdot,\cdot;\mu)$ is coercive and continuous on $\mathbf{V} \times \mathbf{V}$,
%     \item the linear form $f(\cdot;\mu)$ is continuous on $\mathbf{V}$,
% \end{itemize}
Therefore, the problem \eqref{varitional-problem} admits a unique solution.

We introduce a final key assumption essential for the efficiency of the ROM framework: the affine decomposition. This assumption further facilitates the \emph{offline}-\emph{online} paradigm outlined in the Introduction. The bilinear and linear forms introduced are affine if they can be written as:
 % \begin{equation}
 %    \scriptsize
 %     a(u,v;\mu)=\sum_{q=1}^{Q_a} \theta^q_a(\mu)a_q(u,v), \quad
 %     f(v;\mu) = \sum_{q=1}^{Q_f}\theta^q_f(\mu) f_q(v), \quad
 %     l(v;\mu)=\sum_{q=1}^{Q_l} \theta_l^q(\mu)l_q(v),
 %     \label{affine}
 % \end{equation}
\begin{equation}
    \begin{aligned}
    a(u,v;\mu) = \sum_{q=1}^{Q_a} \theta^q_a(\mu)a_q(u&,v),  \quad
    f(v;\mu) = \sum_{q=1}^{Q_f}\theta^q_f(\mu) f_q(v), \\[0.5ex]
    l(v;\mu) &= \sum_{q=1}^{Q_l} \theta_l^q(\mu)l_q(v),
    \end{aligned}
    \label{affine}
\end{equation}
where $a_q$, $f_q$ and $l_q$ are independent on the parameter $\mu \in \mathbf{P}$, while $\theta^q_a$, $\theta^q_f$, $\theta^q_l$ are scalar quantities depending only on the parameter values $\mu \in \mathbf{P}$. The finite dimensional form of equation \eqref{varitional-problem} evaluates $s_{\delta}(\mu) = l(u_{\delta}(\mu);\mu)$, where $u_{\delta}(\mu)$ satisfies
\begin{equation}
      a(u_{\delta}(\mu),v_{\delta};\mu) = f(v_{\delta};\mu) \qquad \forall v_{\delta} \in \mathbf{V}_{\delta}.
      \label{varitional-problem-discr}
\end{equation}
By leveraging the continuity, the coercivity of $a(\cdot,\cdot;\mu)$ and the Galerkin orthogonality, Cea's lemma \cite{monk2003finite} holds:
\begin{equation}
    \lVert u(\mu)- u_{\delta}(\mu) \rVert_{\mathbf{V}} \le \bigg( 1 + C(\mu)%\frac{\gamma(\mu)}{\alpha(\mu)}
    \bigg) \inf_{v_{\delta}\in \mathbf{V}_{\delta}} \lVert u(\mu)-v_{\delta} \rVert_{\mathbf{V}}, \quad \forall v_{\delta} \in \mathbf{V}_{\delta}, \label{eq:Cea1}
\end{equation}
where $C(\mu)$ depends on the coercivity and continuity constant of $a(\cdot,\cdot;\mu)$. As a result, the best approximation error of $u(\mu)$ in the space $\mathbf{V}_{\delta}$ is strongly related to the approximation error $\lVert u(\mu)- u_{\delta}(\mu) \rVert_{\mathbf{V}}$ through the constant $C(\mu)$. The goal is to identify a reduced set of basis functions that can accurately reconstruct the numerical solution $u_{\delta}(\mu)$ through their linear combination. %to find a small number of basis functions whose linear combination exactly represents the numerical solution $u_{\delta}(\mu)$. 

Let $\{ \xi \}_{i=1}^{N}$ be the POD $N$-dimensional set of reduced basis (see Sec. \ref{sec:pod}). Then, the space derived is:
\begin{equation}
    \mathbf{V}_{\text{rb}}=\text{span}\{ \xi_1,\dots,\xi_N \} \subset \mathbf{V}_{\delta},
\end{equation}
where we suppose $N \ll N_{\delta}$. The reduced form of equation \eqref{varitional-problem-discr} evaluates $s_{\text{rb}}(\mu) = f(u_{\text{rb}}(\mu);\mu)$, where $u_{\text{rb}}(\mu)$ satisfies:
  \begin{equation}
      a(u_{\text{rb}}(\mu),v_{\text{rb}};\mu) = f(v_{\text{rb}};\mu) \qquad \forall v_{\text{rb}} \in \mathbf{V}_{\text{rb}}.
      \label{varitional-problem-rb}
  \end{equation}
The reduced solution can be expressed as 
  \begin{equation}
u_{\text{rb}}(\mu)=\sum_{i=1}^{N}(u_{\text{rb}}^{\mu})_i \xi_i,
      \label{varitional-problem-rb2}
  \end{equation}
 where $(u_{\text{rb}}^{\mu})_i$ are the coefficients of the reduced basis approximation. Note that the reduced coefficients depend on the parameter $\mu$, whereas the reduced basis remains parameter-independent. 
 
We now assess the accuracy of the reduced solution. By applying the triangle inequality, we obtain
\begin{equation}
    \lVert u(\mu)- u_{\text{rb}}(\mu) \rVert_{\mathbf{V}} \le \lVert  u(\mu) - u_{\delta}(\mu) \rVert_{\mathbf{V}} + \lVert u_{\delta}(\mu) - u_{\text{rb}}(\mu)  \rVert_{\mathbf{V}}.
\end{equation}
A metric quantifying the distance between $\mathcal{M}_{\delta}$ and $\mathbf{V}_{\text{rb}}$ has been previously introduced in the literature, known as Kolmogorov $N$-width \cite{pinkus2012n}, and is defined as:
\begin{equation}
    d_N(\mathcal{M}_{\delta}) = \inf_{\mathbf{V}_{\text{rb}}} \mathcal{E}(\mathcal{M}_{\delta},\mathbf{V}_{\text{rb}}) = \inf_{\mathbf{V}_{\text{rb}}} \sup_{u_{\delta}\in\mathcal{M}_{\delta}}\inf_{v_{\text{rb}}\in\mathbf{V}_{\text{rb}}} \lVert u_{\delta}-v_{\text{rb}} \rVert_{\mathbf{V}}.
\end{equation}
This quantity serves as a metric for assessing how effectively the features of the system  are captured. In turbulent scenarios, for instance, the flow is decomposed into $N$ finely resolved components to ensure a precise representation. As the Kolmogorov $N$-width increases, the simulation captures more detail, but at the cost of higher computational complexity. Striking a balance between accuracy and computational efficiency is crucial in the modeling process. For elliptic problems, such as diffusion equations, the quantity $d_N(\mathcal{M}_{\delta})$ rapidly decreases as $N$ grows. Consequently, in these cases, a limited number of basis functions is sufficient to accurately approximate the set of high-fidelity solutions $\mathcal{M}_{\delta}$ \cite{pinkus2012n,chen2012certified}. 
Cea's lemma (see eq. \eqref{eq:Cea1}) for $\mathbf{V}_{\text{rb}}$ reads:
\begin{equation}
    \lVert u(\mu)- u_{\text{rb}}(\mu) \rVert_{\mathbf{V}} \le \bigg( 1 + %\frac{\gamma(\mu)}{\alpha(\mu)}
    C(\mu)\bigg) \inf_{v_{\text{rb}}\in \mathbf{V}_{\text{rb}}} \lVert u(\mu)-v_{\text{rb}} \rVert_{\mathbf{V}}, \quad \forall v_{\text{rb}} \in \mathbf{V}_{\text{rb}}.
    \label{cea}
\end{equation}
This result highlights that the ROM accuracy is determined by the specific properties of the problem under consideration. Eq. \eqref{varitional-problem-rb} can be formulated in matrix form: evaluate $s_{\text{rb}}(\mu) = f(u_{\text{rb}}(\mu);\mu)$, where $u_{\text{rb}}(\mu)$ satisfies
  \begin{equation}
      a(u_{\text{rb}}(\mu),\xi_j;\mu) = f(\xi_j;\mu) \qquad  1 \le j \le N.
      %\label{varitional-problem-rb}
  \end{equation}
By using eq. \eqref{varitional-problem-rb2} %$u_{\text{rb}}(\mu)=\sum_{i=1}^{N}(u_{\text{rb}}^{\mu})_i \xi_i$,  
it becomes: evaluate $s_{\text{rb}}(\mu) = {(u^{\mu}_{\text{rb}})}^T f^{\mu}_{\text{rb}}$ such that
\begin{equation}
A^{\mu}_{\text{rb}}u^{\mu}_{\text{rb}}=f^{\mu}_{\text{rb}},
\end{equation}
where
\begin{equation}
    (A^{\mu}_{\text{rb}})_{j,i} =  a(\xi_i,\xi_j;\mu) \qquad \text{and} \qquad (f^{\mu}_{\text{rb}})_j = f(\xi_j;\mu), \qquad 1 \le i,j \le N.
\end{equation}
To conclude, projection-based ROMs offer significant advantages, primarily due to their direct derivation from the original governing equations, ensuring a strong connection to the underlying physics of the problem. We expect that this approach preserves essential physical properties, provides a coherent numerical solution throughout the simulation, and remains effective even in data-scarce scenarios \cite{amsallem2012stabilization,balajewicz2012stabilization,karatzas2020projection}.

\subsection{Data-driven models}
\label{sec:data-driven}

Unlike projection-based methods, data-driven ROMs do not rely on the governing equations or any prior physical knowledge of the system. Typically, these methods achieve speed-ups several orders of magnitude greater than those of intrusive ROMs \cite{aromabook}.
As a consequence, these methods are particularly beneficial for real-time applications requiring rapid decision-making, such as control systems and optimization. Conversely, the accuracy of these methods relies heavily on the quality and quantity of the available data, especially their ability to capture the underlying dynamics. Additionally, non-intrusive ROMs lack a well-established error estimation theory. In data-driven ROMs, techniques such as regression models, Neural Networks (NNs) and Gaussian processes are widely employed to extract patterns and relationships from data \cite{milano2002neural,chen2021physics,hesthaven2018non,dar2023artificial}. These methods enable the identification of key features and correlations that may be difficult to capture using conventional modeling approaches, particularly in systems characterized by strong nonlinearities \cite{papapicco2021neural,fresca2022pod,boukraichi2023parametrized}. While the integration of data, machine learning techniques, and dimensionality reduction facilitates the development of accurate and flexible models, several challenges remain open. These include the necessity of representative datasets to ensure model reliability and the mitigation of overfitting \cite{daniel2020model}, particularly in scenarios with limited data availability.

To further explore this framework, we now introduce data-driven ROMs in the context of our study.
Unlike intrusive techniques, the coefficients $u_{\text{rb}}^{\mu}$ are not obtained by solving a differential-algebraic system derived from the original governing equations. Instead, they are determined through a mapping in the parameter space.
Specifically, given the solution manifold (or a database that may also incorporate experimental data), the reduced coefficients $u_{\text{rb}}^{\mu}$ are computed by projecting the high-fidelity snapshots onto the reduced space. The set $\{ (\mu, u_{\text{rb}}^{\mu}) : \mu \in \mathbf{P}_h \}$ is approximated through various techniques, such as NNs, radial basis functions, and Gaussian process regression, among others. In the following, we focus on NNs, as they are employed in Sec.~\ref{sec:aorta}.

\subsection{Neural networks}
\label{sec:neural-network}

A neural network is a computational model designed to learn patterns and relationships from data. One of its key advantages is its capability as a \emph{universal approximator} \cite{cybenko1989approximation}, enabling it to capture complex non-linear dependencies through an iterative learning process. This property is particularly beneficial in ROM when dealing with a wide range of physical and geometrical configurations, as often encountered in real-world applications. However, significant challenges remain, including the selection of an optimal network architecture, the computational cost associated with repeated training, and the substantial amount of data required for effective learning. A standard NN is composed of a set of neurons interconnected through directional, weighted synaptic connections. In this structure, neurons function as nodes, while synapses represent edges within an oriented graph. Each neuron $j$ is characterized by three functions:
\begin{itemize}
    \item the propagation function $u_j$ defined as: 
    $$u_j=\sum_{k=1}^m w_{s_k,j}y_{s_k} + b_j,$$ where $b_j$ represents the bias, $y_{s_k}$ denotes the output of the sending neuron $k$,  $w_{s_k,j}$ are the weights and $m$ is the number of sending neurons connected to the neuron $j$. 

    \item the activation function $a_j$ given by: 
    $$ a_j=f_{\text{act}}\left(\sum_{k=1}^m w_{s_k,j}y_{s_k}+b_j\right).$$
    Common choices for $f_{\text{act}}$ include sigmoid function, hyperbolic tangent, ReLU and SoftMax \cite{sharma2017activation}. 
    \item the output function $y_j$ identified as:
    \begin{equation*}
        y_j=f_{\text{out}}(a_j). %=a_j.
    \end{equation*}
    In many cases, the output function is simply the identity function, leading to $y_j = a_j$.
\end{itemize}
Feedforward NNs, in which neurons are organized into sequential layers, are commonly employed for interpolation purposes \cite{rosenblatt1958perceptron,fine2006feedforward}. These networks consist of an input layer, an output layer and one or more hidden layers in between. During the training process, the network weights are iteratively adjusted to minimize the discrepancy between the predicted outputs and the actual values from the training dataset. This discrepancy is quantified using a loss function, denoted as $\mathcal{L}$. For interpolation problems, a common choice for $\mathcal{L}$ is the mean squared error. The backpropagation algorithm \cite{rojas1996backpropagation,rumelhart1986learning} is employed to update the weights, leveraging the gradient of $\mathcal{L}$ with respect to the network parameters:
\begin{equation}
\begin{split}
    & \frac{\partial \mathcal{L}}{ \partial w_{s_k,j}^l}=\frac{\partial \mathcal{L}}{\partial a_j^l}\frac{\partial a_j^l}{\partial u_j^l}\frac{\partial u_j^l}{\partial w_{s_k,j}^l} , \\
    & \frac{\partial \mathcal{L}}{ \partial b_j^l}=
   \frac{\partial \mathcal{L}}{\partial a_j^l}\frac{\partial a_j^l}{\partial u_j^l} \frac{\partial u_j^l}{\partial b_j^l}.
\end{split}
\end{equation}
This algorithm enables the computation of the error gradient with respect to the weights for a given input by propagating the error backward through the network. Specifically, the forward pass calculates the output layer values starting from the input layer and evaluates the loss function. Subsequently, backpropagation executes a backward pass to determine the gradient of the loss function and update the model parameters accordingly:
\begin{equation}
\begin{split}
    & \bm{w} = \bm{w} - \eta \frac{\partial \mathcal{L}}{\partial \bm{w}}, \\
    & \bm{b} = \bm{b} - \eta \frac{\partial \mathcal{L}}{\partial \bm{b}},
\end{split}
\end{equation}
where $\bm{w}$ and $\bm{b}$ represent the weight and bias matrices, respectively, while $\eta$ denotes the learning rate. Proper tuning of hyperparameters -- including the learning rate, the activation function, and the number of hidden layers and neurons -- is crucial to optimize the performance of the network \cite{goodfellow2016deep,kriesel2007brief,calin2020deep}.

%\textcolor{black}{Aggiungiamo qualcosa sugli hybrid models?}

\section{A patient-specific application}
\label{sec:aorta}
\textcolor{black}{In this section, we present a patient-specific cardiovascular application that combines Galerkin projection with data-driven techniques within a hybrid ROM framework. Hybrid approaches are particularly advantageous when the physical model is either partially known or computationally expensive to evaluate. In these cases, the dominant dynamics are captured using reduced-order equations, while specific components, such as nonlinear source terms, stabilization mechanisms, effective parameters or closure models, are learned from data. These data-driven elements can be approximated using interpolation, regression and machine learning techniques, resulting in reduced models that are both efficient and generalizable. This strategy preserves the interpretability of physics-based models while leveraging data to account for unknown or complex features. Moreover, it enables the extension of traditional equation-based ROM frameworks to more realistic and demanding scenarios, overcoming limitations associated with purely physics-based or black-box machine learning methods.}
\subsection{Problem formulation}
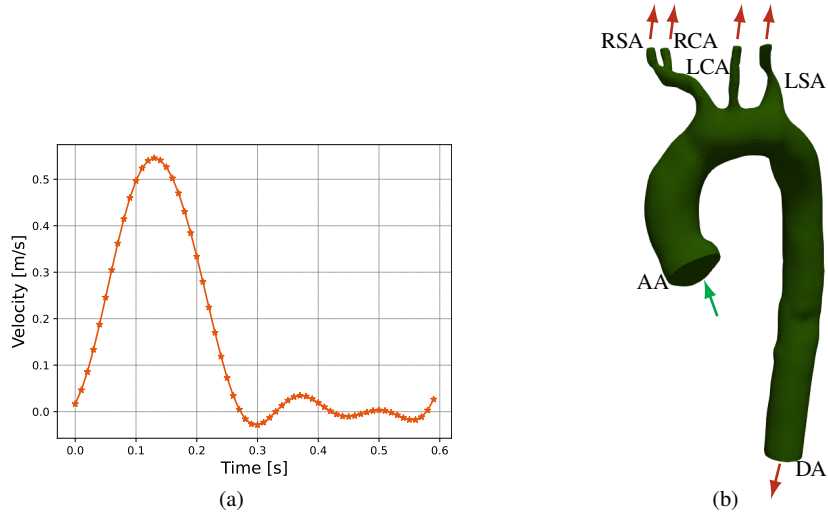
\begin{figure}
	\centering
    \subfloat[][\label{fig:BC_aorta}]{\includegraphics[width=.5\textwidth]{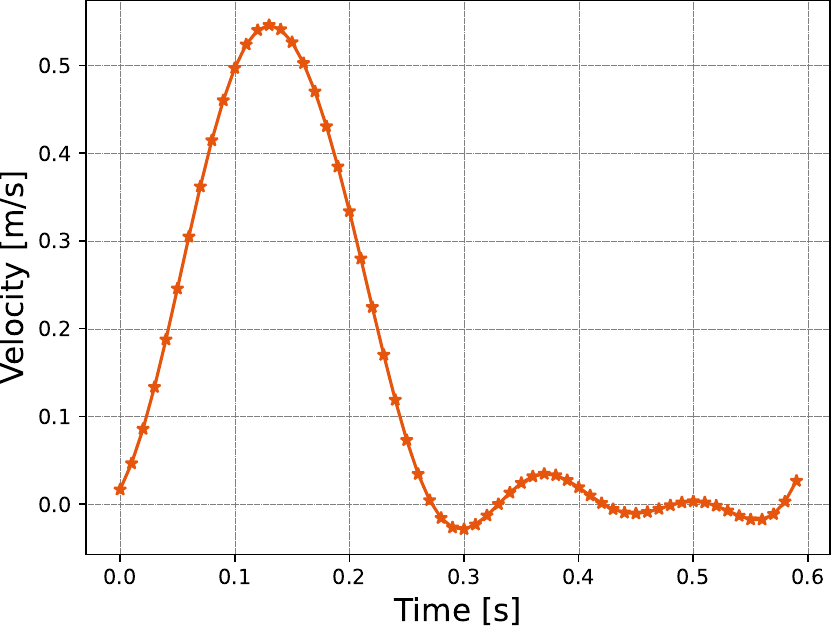}}\hspace{2cm}
    \subfloat[][\label{fig:aorta_mesh}]{
    %\subfloat[][\label{fig:BC_aorta}]{\includegraphics[width=.4\textwidth]{img/BC_velocity_aorta.pdf}}\hspace{2cm}
 	%\subfloat[][\label{fig:aorta_mesh}]{%\includegraphics[width=.13\textwidth]{img/aorta_mesh.png}
    \begin{overpic}[width=0.25\textwidth]{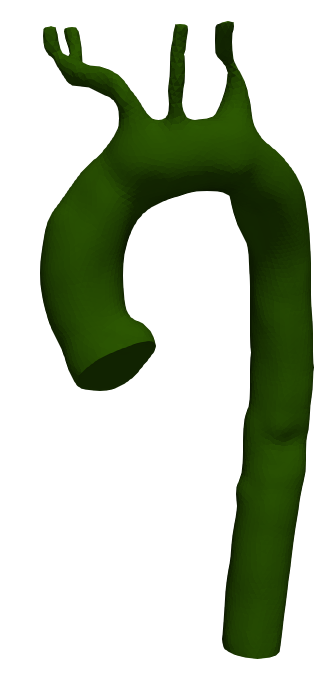}

    \put(4,42){\scalebox{1}{AA}}
    \put(39.5,0.2){\scalebox{1}{DA}}
    \put(-4,96){\scalebox{1}{RSA}}
    \put(12.2,96){\scalebox{1}{RCA}}
    \put(15,90.2){\scalebox{1}{LCA}}
    \put(37,87){\scalebox{1}{LSA}}
    
    \linethickness{1pt}\put(22,36){\color{Green}\vector(-0.08,0.22){3}}
    
    \linethickness{1pt}\put(36,3){\color{BrickRed}\vector(-0.05,-0.20){2}}

     \linethickness{1pt}\put(11.5,98){\color{BrickRed}\vector(0.01,0.08){1}}

     \linethickness{1pt}\put(7,98){\color{BrickRed}\vector(0.01,0.08){1}}

     \linethickness{1pt}\put(26.5,98){\color{BrickRed}\vector(0.01,0.08){1}}

     \linethickness{1pt}\put(33,98){\color{BrickRed}\vector(0.01,0.08){1}}
    \end{overpic}
    %}\\
    }
	\caption{Time evolution of the boundary condition for the velocity $\bm u_D(t)$ (a) and sketch of the computational domain $\Omega$ (b).}
	\label{fig:BC_aorta_mesh}
\end{figure}

We analyze the blood flow dynamics within a patient-specific aortic arch. The domain is shown in Figure~\ref{fig:aorta_mesh}. The boundary consists of an inlet section, marked by a green arrow, representing the Ascending Aorta (AA), and five outlet sections -- Right Subclavian Artery (RSA), Right Common Carotid Artery (RCA), Left Common Carotid Artery (LCCA), Left Subclavian Artery (LSA), and Descending Aorta (DA) -- all indicated by red arrows. 

The blood is modeled as an incompressible Newtonian fluid, described by the unsteady Navier-Stokes equations within the spatial domain $\Omega \in \mathbb{R}^3$ over the time interval $(t_0, T]$:
   \begin{empheq}[left=\empheqlbrace]{alignat=4}
       \frac{\partial \bm{u}}{\partial t} + \nabla \cdot \left(\bm{u} \otimes \bm{u}\right) - \nabla \cdot (\nu \nabla \bm{u}) +\nabla p & =  0 & \quad \mbox{ on } \Omega \times (t_0, T], \label{eq:NS1} \\
      \nabla \cdot \bm{u} & =  0 & \quad \mbox{ on } \Omega \times (t_0, T], \label{eq:NS2}
   \end{empheq}
where $\bm{u}=\bm{u}(\bm{x}, t)$ is the velocity vector, $p=p(\bm{x}, t)$ is the kinematic pressure, and $\nu$ is the kinematic viscosity.

Time-dependent Dirichlet conditions are imposed at the inlet (see the velocity profile in Figure \ref{fig:BC_aorta}), no-slip constraints are enforced along the walls, and homogeneous Neumann conditions are applied at the outlet sections:
\begin{equation}
\begin{cases}
    \bm u = \bm u_D(t) &\quad \text{in $\Gamma_i\times (t_0, T)$}, \\
    \bm u  = \bm 0 &\quad \text{in $\Gamma_w\times (t_0, T)$}, \\
    (2\mu \nabla \bm u - p I)\bm n = \bm 0 &\quad \text{in $\Gamma_o\times (t_0, T)$}.
\end{cases}
\label{BC0}
\end{equation}

% For the velocity field a nonhomogeneous Dirichlet boundary condition is enforced on the inlet section:
% \begin{equation}
%     \label{eq:non-hom_dirich}
%     \bm{u} = \bm{u}_D \quad \mbox{ on } \Gamma_{i} \times (t_0, T],
% \end{equation}
% homogeneous Neumann boundary conditions are applied on the outlet sections:
% \begin{equation}
%     \nabla \bm{u} \cdot \bm{n} = 0 \quad \mbox{ on } \Gamma_{o} \times (t_0, T],
% \end{equation}
% and finally a no-slip boundary condition is employed on the wall
% \begin{equation}
%     \bm{u} = 0 \quad \mbox{ on } \Gamma_{w} \times (t_0, T].
% \end{equation}

\begin{figure}
\begin{center}
\begin{circuitikz}[american,scale=0.65]
%\draw (1,0) to[R=$R_{p,k}$] (2,0);
\draw[] (-0.5,0.5) -- (-0.50,-0.5);
%\draw [-to](0,0) -- (0.5,0);
\draw [-{Stealth[scale=1.5]}] (-0.5,0) -- (0.3,0);
\draw[color=violet] (0,0) to[R, l^=$R_{p}$] (3,0);
\draw[color=lime!70!black] (3,0) to[C, l=$C$,-] (3,-2.5);
\draw[color=teal!50!black] (3,0) to[R, l^=$R_{d}$] (6,0);
%\ctikzset{tripoles/mos style/arrows}
\draw[] (0,0.4) node[] {$Q$};
\draw (-0.8,0) node[] {$p$};
\draw (3,0.2) node[] {$p_{p}$};
\draw (3,-2.7) node[] {$p_{d}$};
\draw (6.3,0) node[] {$p_{d}$};
\draw[] (-0.45,-0.8) node[] {$\Gamma_o$};
\end{circuitikz}
\end{center}
\caption{Sketch of the three-element Windkessel model.}
\label{fig:windkessel}
\end{figure}
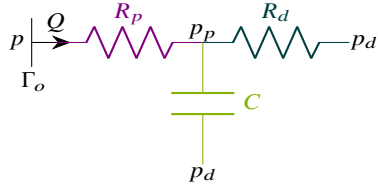

% On the other hand, for the pressure field, a homogeneous Neumann boundary condition is enforced on the inlet and on the wall sections:
% \begin{equation}
%     \nabla p \cdot \bm{n} = 0 \quad \mbox{ on } (\Gamma_{w} \cup \Gamma_{i}) \times (t_0, T].
% \end{equation}
In order to enforce realistic outflow boundary conditions for the pressure, at each outlet of the domain, a three-element Windkessel model is employed
\cite{westerhof2009arterial}. This model is composed of a compliance $C$ and two resistances, namely the proximal resistance $R_p$ and the distal resistance $R_d$. On a specific outlet section, the downstream pressure $p$ is given by the following differential algebraic equations system: 
\begin{equation}
   \left \{
   \begin{alignedat}{3}
    C \frac{dp_p}{dt} + \frac{p_p - p_d}{R_d}  & = Q & \quad & \mbox{ on } \Gamma_o,\\
    p - p_p & = R_p Q & \quad & \mbox{ on } \Gamma_o,
   \end{alignedat}
   \right .
   \label{eq:Windkessel-FOM}
\end{equation}
where $p_p$ is the proximal pressure, $p_d$ is the distal pressure (assumed to be, as common in literature \cite{nichols2022mcdonald}, $p_d=0$, as it serves as a reference value) and $Q$ is the flow rate through the outlet section. A scheme of the Windkessel model is shown in Figure \ref{fig:windkessel} and the Windkessel coefficients for each outlet are shown in Table \ref{table_windkessel_aorta}.
Both the inflow boundary condition (Figure~\ref{fig:BC_aorta}) and the Windkessel coefficients are derived from RHC and ECHO tests in \cite{Girfoglio2020}. The kinematic viscosity is $\nu=3.7\cdot 10^{-6}$ m$^2$/s.

\begin{table}[htb!] 
\caption{Parameter values of the Windkessel model.}
\centering
\renewcommand{\arraystretch}{1.5}
\begin{tabular}{|c|c|c|c|}
\hline
\rowcolor{gray!60} 
 Outlet & $R_p$ [m$^{-1}$s$^{-1}$] & $R_d$ [m$^{-1}$s$^{-1}$] & $C$ [ms$^{2}$] \\
\hline
Right subclavian artery (RSA) & 1.84$\cdot 10^{8}$  & 3.11$\cdot 10^{9}$ & 3.26$\cdot 10^{-5}$ \\
\hline\rowcolor{gray!20} 
Right common carotid artery (RCA) & 1.23$\cdot 10^{8}$ & 2.07$\cdot 10^{9}$  & 5.16$\cdot 10^{-10}$ \\
\hline
Left common carotid artery (LCA) & 1.78$\cdot 10^{8}$ & 3.01$\cdot 10^{9}$ & 3.52$\cdot 10^{-10}$ \\
\hline\rowcolor{gray!20} 
Left subclavian artery (LSA) & 7.09$\cdot 10^{7}$  & 1.19$\cdot 10^{9}$  & 9.35$\cdot 10^{-10}$ \\
\hline
Descending aorta (DA) & 7.8$\cdot 10^{6}$ & 1.31$\cdot 10^{8}$ & 7.72$\cdot 10^{-9}$ \\
\hline
\end{tabular}
\label{table_windkessel_aorta}
\end{table}

For the space discretization, we employ the finite volume method and the time derivative is approximated with a first-order Euler scheme. In order to treat the velocity-pressure coupling, the Pressure Implicit with Splitting of Operators (PISO) algorithm \cite{issa1986solution} is adopted. The Euler scheme is also used to discretize the Windkessel system \eqref{eq:Windkessel-FOM}.

\subsection{The reduced framework}
\label{sec:rom-aorta}
The ROM adopted is a hybrid approach, which combines the POD-Galerkin procedure, the lifting function method \cite{fick2018stabilized,graham1999optimal} and NNs. The time is the only parameter involved. \textcolor{black}{While only the temporal parameter is considered here, the same procedure can be extended to include geometric variations, with snapshots collected for different shapes. This extension, however, typically increases the offline computational cost and may slightly deteriorate the ROM accuracy. A sensitivity analysis on the number of snapshots could help in mitigating this effect.} %\textcolor{black}{In future developments, the framework can be extended to include geometric parameters as well. In this case, once a geometric parameter to be varied is selected, one can collect snapshots at different values of this parameter and apply the same methodology described in this chapter, extending the dependence of the solution from time to both time and geometry. Based on previous experience, a slight deterioration of the ROM accuracy is expected when geometric variations are introduced. A possible direction would be to perform a sensitivity analysis with respect to the number of snapshots, in order to optimize the computational cost of the offline phase, which typically increases in the presence of geometric parameters.} 

The key assumption is that the solution can be approximated as a linear combination of spatially dependent basis functions for velocity and pressure, respectively $\bm{\phi}(\bm{x})$ and $\psi(\bm{x})$, and time-dependent coefficients, respectively $a_i(t)$ and $b_i(t)$: \begin{equation}
    \bm u \approx \bm u_{\text{rb}} = \sum_{i=1}^{N_{\bm u}} a_i(t) \bm\phi_i(\bm x),
    \label{u_rb}
\end{equation}
\begin{equation}
    p \approx p_{\text{rb}} = \sum_{i=1}^{N_p} b_i(t) \psi_i(\bm x),
    \label{p_rb}
\end{equation}
where $N_p$ and $N_u$ are the dimensions of the reduced basis space for pressure and velocity, respectively. This assumption enables the separation of computations into a computationally intensive offline phase, executed only once, and a computationally efficient online phase, performed for each new evaluation. During the offline stage, the FOM solutions are computed for a discrete set of time instances, followed by the computation of lifting functions to homogenize the snapshots. The POD is then applied, and a Galerkin projection is performed to derive the reduced-order model. Additionally, a feedforward neural network is trained to generalize the outflow pressure dynamics across different time values. In the online phase, for a new set of time instances, the corresponding modal coefficients are determined by solving the reduced dynamical system obtained in the offline phase. Finally, the ROM solution is reconstructed by combining these coefficients with the reduced basis functions.

In order to compute the lifting functions $\bm\chi^{\bm u}$ and $\bm\chi^{p}$ we solve a potential flow problem.
% \begin{equation}
%     \begin{cases}
%     \begin{aligned}
%         %\nabla \cdot \bm u = 0 & \quad \mbox{ on } \Omega, \\
%         \Delta \bm\chi^{\bm u} &= 0  \quad \mbox{ on } \Omega,\\ %\bm u = \nabla \bm\chi^{\bm u} & \quad \mbox{ on } \Omega, \\
%         \bm\chi^{\bm u} &= 1  \quad \mbox{ on } \Gamma_i,\\
%         \bm\chi^{\bm u} &= 0  \quad \mbox{ on } \Gamma_w,\\
%         \nabla \bm\chi^{\bm u} \cdot \bm n &= 1  \quad \mbox{ on } \Gamma_o,
%     \end{aligned}
%     \end{cases}\qquad \mbox{and} \qquad 
%     \begin{cases}
%     \begin{aligned}
%         \Delta \bm\chi^p + \nabla\cdot ( \nabla \cdot (\bm\chi^{\bm u} \otimes \bm\chi^{\bm u})) &= 0   \quad \mbox{ on } \Omega, \\  
%         \bm\chi^{p} &= 1  \quad \mbox{ on } \Gamma_o,\\
%         \nabla \bm\chi^{p} \cdot \bm n &= 0  \quad \mbox{ on } \partial \Omega \setminus \Gamma_o.
%     \end{aligned}
%     \end{cases}
%     \label{potentialflow}
% \end{equation}
Then each snapshot is modified as follows:
\begin{equation}
\bm u' (\bm x, t) = \bm u (\bm x, t) - 
{\bm u}_D(t) \bm\chi^{\bm u}(\bm x),
\label{shiftu}
%\sum_{b=1}^{N^{\bm u}_{BC}} g^{\bm u}_b(t) \bm\chi^{\bm u}_b(\bm x),
\end{equation}
\begin{equation}
 p' (\bm x, t) =  p (\bm x, t) - 
 {p}_D(t) \bm\chi^{p}(\bm x),
 %\sum_{b=1}^{N^{p}_{BC}} g^{p}_b(t) \bm\chi^{p}_b(\bm x),
 \label{schiftp}
\end{equation}
where ${\bm u}_D(t)$ is the prescribed boundary condition for the velocity on $\Gamma_i$ 
and ${p}_D(t)$ is the pressure outflow on $\Gamma_{o}$ computed discretizing \eqref{eq:Windkessel-FOM} and interpolating these values with a NN. This transformation is essential to handle non-homogeneous boundary conditions at the ROM level.

The $L^2$ orthogonal projection of equation \eqref{eq:NS1} %
onto $\text{span}\{\bm\phi_1,\dots,\bm\phi_{N_{\bm u}}\}$ results in \cite{akhtar2009stability,bergmann2009enablers,lorenzi2016pod}:
\begin{equation}
    \big(\bm \phi_i, \partial_t \bm{u} + \nabla \cdot (\bm{u} \otimes \bm{u}) - \nabla \cdot (\nu \nabla \bm{u})+\nabla p \big)_{L^2(\Omega)} = 0 \quad \text{for }i = 1, \dots, N_{\bm u}.
    \label{proj_l2}
\end{equation}
By substituting equations \eqref{u_rb} and \eqref{p_rb} in equation \eqref{proj_l2} and given the orthonormality of the reduced basis, we get the following dynamical system: %for the reduced coefficients:
\begin{equation}
    \dot{\bm a} = \nu \bm B \bm a - \bm a^T \bm C \bm a - \bm K \bm b,
    \label{galerkin_u}
\end{equation}
where $\bm a = \{ a_i(t)\}_{i=1}^{N_u}$, $\bm b = \{ b_i(t)\}_{i=1}^{N_p}$, and %are the unknown ROM coefficients for velocity and pressure and:
\begin{equation}
    B_{ij} = \big( \bm\phi_i, \Delta \bm\phi_j \big)_{L^2(\Omega)},
\end{equation}
\begin{equation}
    C_{ijk} = \big( \bm\phi_i, \nabla \cdot (\bm\phi_j \otimes \bm\phi_k) \big)_{L^2(\Omega)},
\end{equation}
\begin{equation}
    K_{ij} = \big( \bm\phi_i, \nabla \psi_j \big)_{L^2(\Omega)}.
\end{equation}
The continuity equation \eqref{eq:NS2} is projected as well onto $\text{span}\{\psi_1, \dots, \psi_{N_p}\}$, providing: %
\begin{equation}
    \big( \psi_i,\nabla \cdot \bm u \big)_{L^2(\Omega)} = 0 \quad \text{for }i = 1, \dots, N_{p}.
    \label{proj_p_continuity}
\end{equation}
By substituting  equation \eqref{u_rb} in equation \eqref{proj_p_continuity}, we get:
\begin{equation}
    \bm P {\bm a} = 0,
    \label{galerkin_p_continuity}
\end{equation}
where
\begin{equation}
    P_{ij} = \big( \bm\psi_i, \nabla \cdot\bm\phi_j \big)_{L^2(\Omega)}.
\end{equation}
To ensure the inf-sup condition at the reduced level, two main strategies are proposed in the literature: supremizer enrichment and the Pressure Poisson Equation (PPE) approach \cite{ballarin2015supremizer}. In this work, we choose the supremizer approach, as it has demonstrated greater effectiveness for this specific case. Both methods provide similar performance in terms of solution accuracy; however, the supremizer approach offers advantages in terms of numerical stability and robustness \cite{stabile2018finite}, particularly in scenarios with moderate Reynolds numbers. Therefore, in this work, we adopt the supremizer approach.

In this study, the inlet velocity is expressed as a continuous function for each time value, providing readily accessible data. In contrast, the outlet pressure is not available continuously; it can only be obtained at discrete time intervals from the numerical solution of the Windkessel model. Consequently, the NN is introduced to establish a functional representation of the outlet pressure. The NN is designed to accept time as an input variable, allowing it to effectively predict the outlet pressure, once it has been properly trained. 

This reduced approach introduces key innovations, notably the use of the lifting function for pressure, which has previously been applied only to velocity and never in combination with the Windkessel model and non-homogeneous pressure. Furthermore, the integration of neural networks marks a significant advancement, enabling the efficient evaluation of pressure values at arbitrary time instants.

\subsection{Numerical results}

The system reaches a pseudo steady-state after 8 cardiac cycles (4.8 s), indicating transient effects have subsided. We focus our analysis on the cycle [5.4,6] s. To train the ROM, we collect 100 snapshots over this interval, using a reduced time step of $\Delta t_r=0.01$. Therefore, the selected time instances are:

\begin{equation}
    \left\{5.4,5.41,5.42,\dots,5.98,5.99,6\right\} s.
    \label{time_set}
\end{equation}

\begin{figure}
	\centering
    \subfloat[][\label{fig:eigup}]{\includegraphics[width=.4\textwidth]{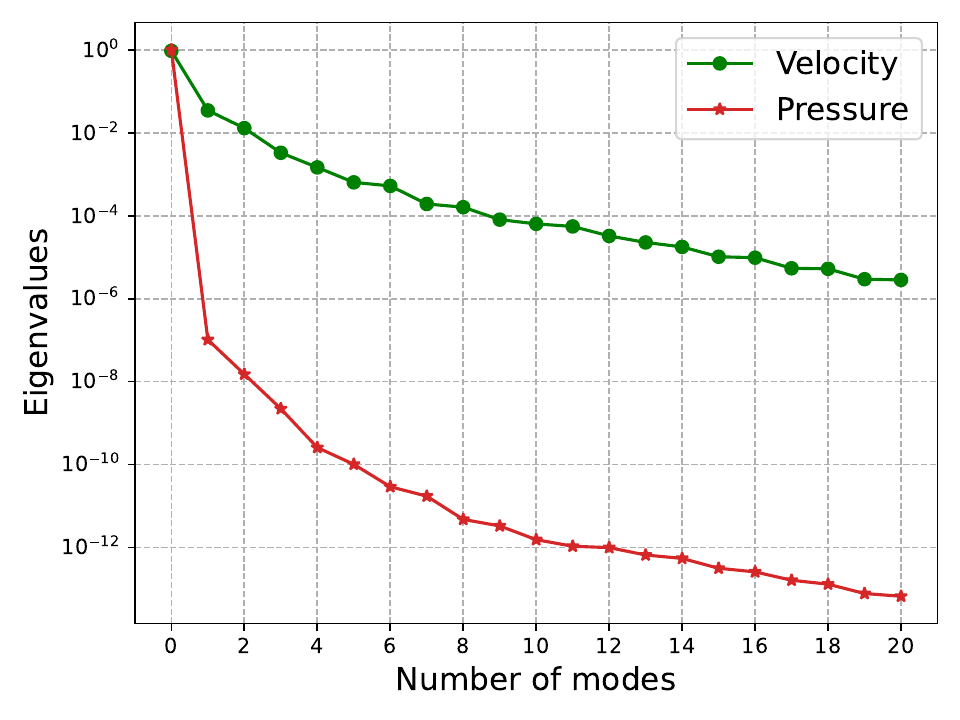}}
 	\hspace{2ex}\subfloat[][\label{fig:cumeigup}]{\includegraphics[width=.4\textwidth]{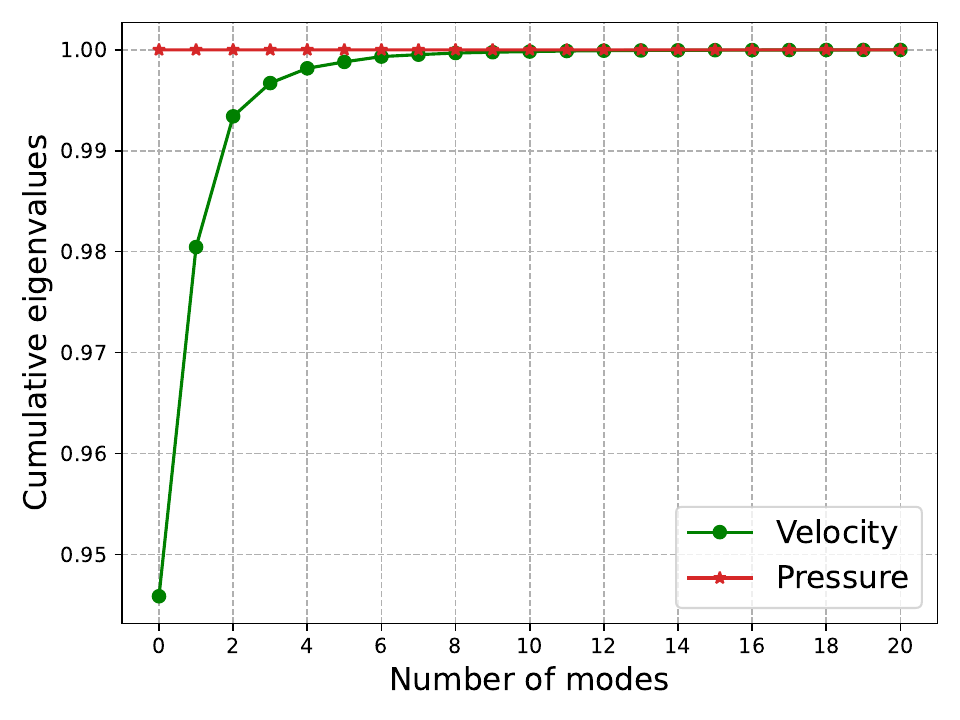}}\\
	\caption{Decay of the eigenvalues (a) and cumulative energy (b) for the velocity and pressure fields.}
	\label{fig:eigup_cumeig}
\end{figure}

\textcolor{black}{In Figure \ref{fig:eigup}, we show the decay of the eigenvalues as the number of modes increases.} To accurately represent $99.99\%$ of the cumulative energy of the eigenvalues \textcolor{black}{(see Figure \ref{fig:cumeigup})}, 12 modes are necessary for the velocity field, while only a single mode suffices for the pressure. 

% \begin{figure}[!h]
%     \centering
%     \includegraphics[width=0.5\linewidth]{img/eig_up.pdf}
%     \caption{Decay of the eigenvalues for velocity and pressure.}
%     \label{fig:eigup_cumeig}
% \end{figure}

\begin{figure}
	\centering
    \subfloat[]{\includegraphics[width=.4\textwidth]{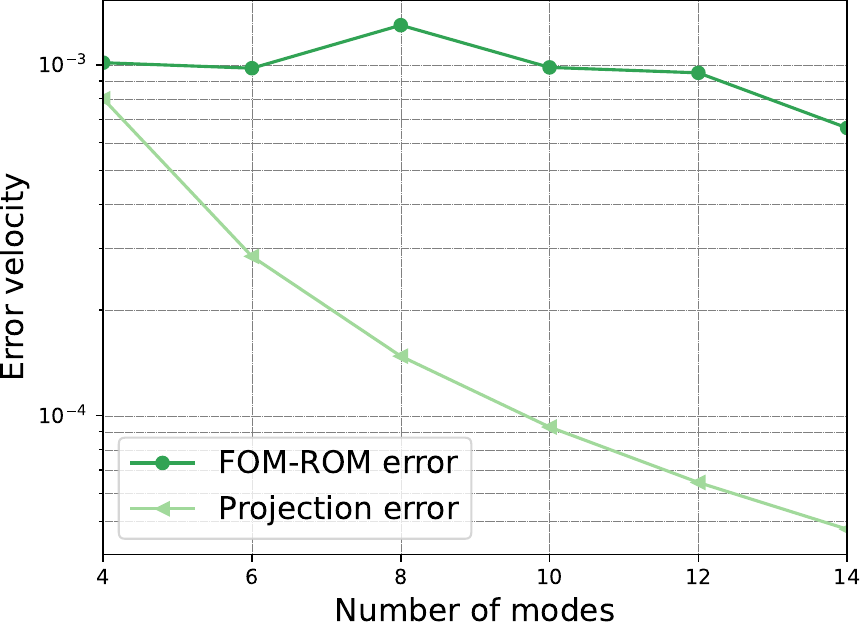}}
 	\subfloat[]{\includegraphics[width=.4\textwidth]{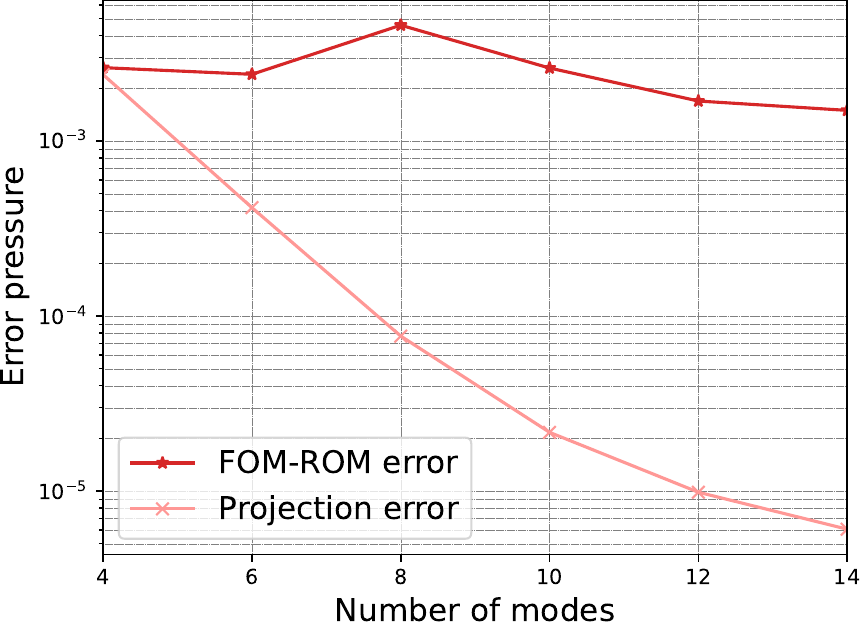}}\\
	\caption{Time-averaged reconstruction and projection error for velocity (a) and pressure (b) as the number of modes $N$ increases.}
    %\includegraphics[width=.6\textwidth]{img/err_aorta_N.pdf}
	%\caption{Time-averaged reconstruction error for velocity and pressure as the number of modes $N$ increases.}
	\label{errup_N_aorta}
\end{figure}

%\textcolor{black}{As depicted in Figure \ref{errup_N_aorta}, the time-averaged reconstruction error slightly decreases between N=4 and N=6, where it reaches its minimum value. For larger numbers of modes, the error increases again, likely due to the inclusion of low-energy modes that amplify numerical noise. Hence, N=6 modes are selected as the optimal number for the reduced basis.}
\textcolor{black}{As depicted in Figure \ref{errup_N_aorta}, the time-averaged reconstruction error decreases up to $N=6$, after which it rises again, surpassing the initial value, before decreasing once more. Therefore, $N=6$ is selected as the optimal number of modes. We clarify that $N$ now refers to each component separately.}

\begin{figure}
	\centering
    \subfloat[][\label{fig:erru_aorta_proj}]{\includegraphics[width=.4\textwidth]{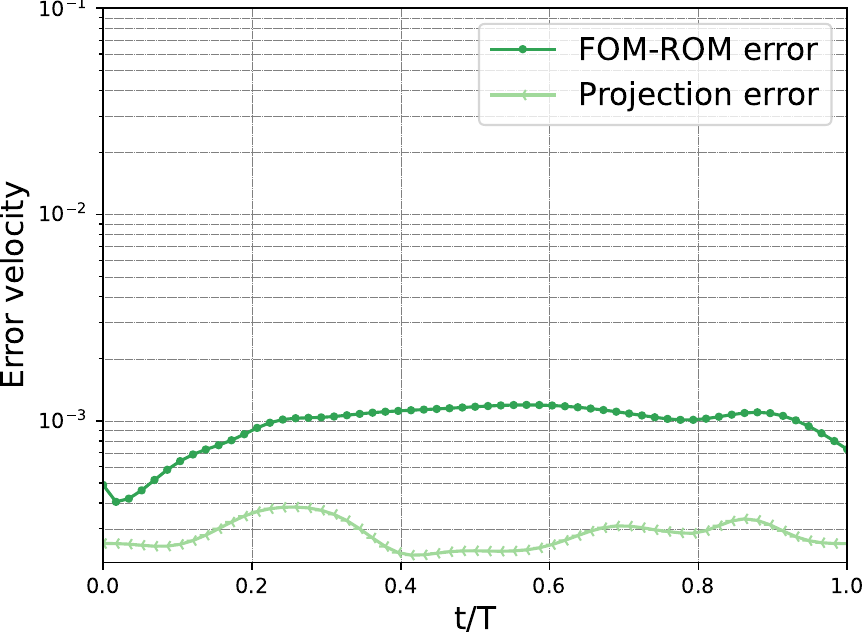}}\hspace{2ex}
 	\subfloat[][\label{fig:error_dtp_proj}]{\includegraphics[width=.4\textwidth]{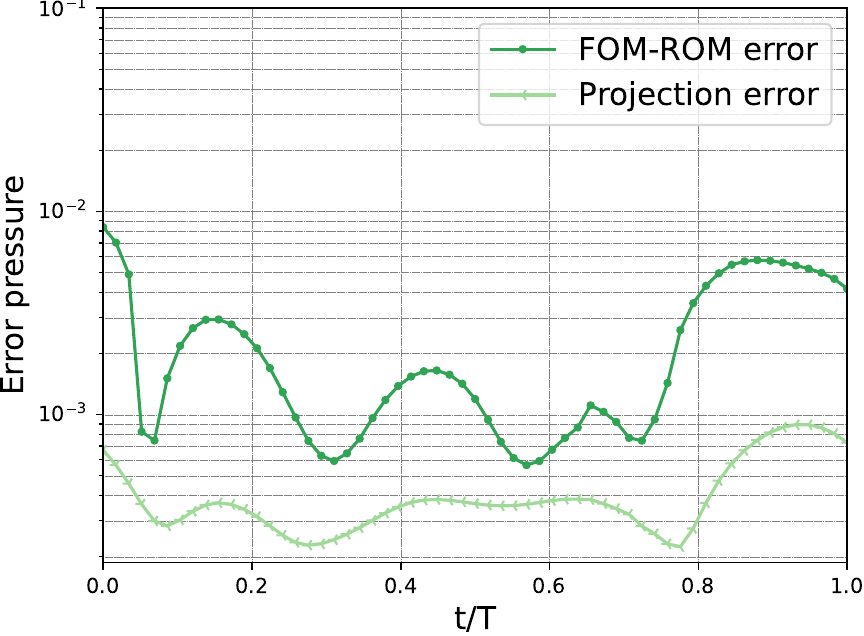}}\\
	\caption{Comparison between reconstruction error and projection error for velocity and pressure with $N = 6$. The supremizer approach is adopted and $N$ is the same for pressure, velocity and supremizers.}
	\label{errup_aorta_proj}
\end{figure}

As an initial numerical experiment, we do not introduce new test points during the online phase, meaning the validation set is identical to the training set \eqref{time_set}.
Figure~\ref{errup_aorta_proj} reports the comparison between the reconstruction error and the projection error.
For both velocity and pressure, the reconstruction error closely aligns with the projection error, confirming the effectiveness of our ROM approach in accurately capturing the flow dynamics. Specifically, the velocity reconstruction error remains on the order of $10^{-3}$ throughout the entire time window, while the pressure error varies between $10^{-3}$  and $10^{-2}$. Note that we calculate the absolute reconstruction error for both velocity and pressure. However, qualitative comparisons are also presented to illustrate that the absolute error serves as an effective metric to assess the quality of the ROM reconstruction.

\begin{figure}
	\centering
    \subfloat[][\label{fig:erru_aorta_noliftP}]{\includegraphics[width=.4\textwidth]{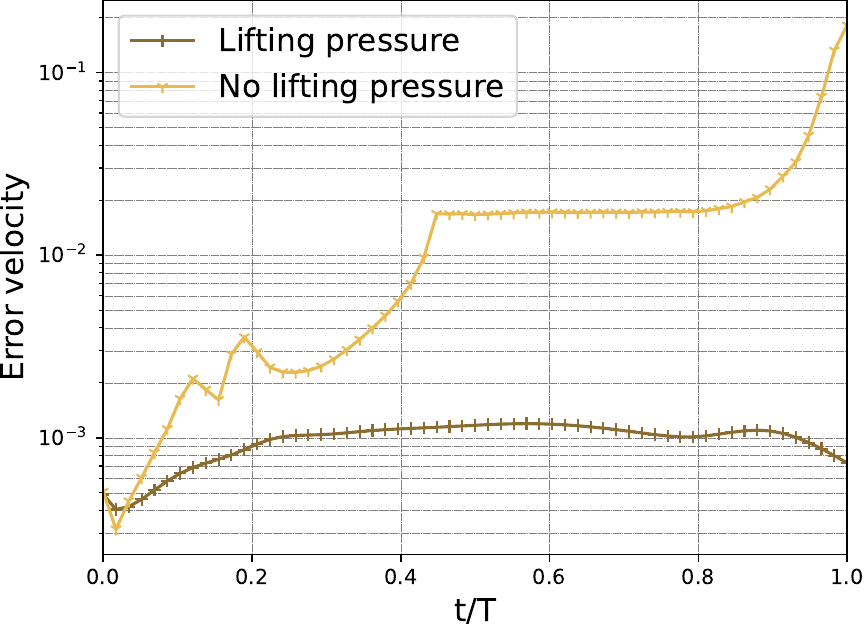}}
 	\hspace{2ex}\subfloat[][\label{fig:error_dtp_noliftP}]{\includegraphics[width=.4\textwidth]{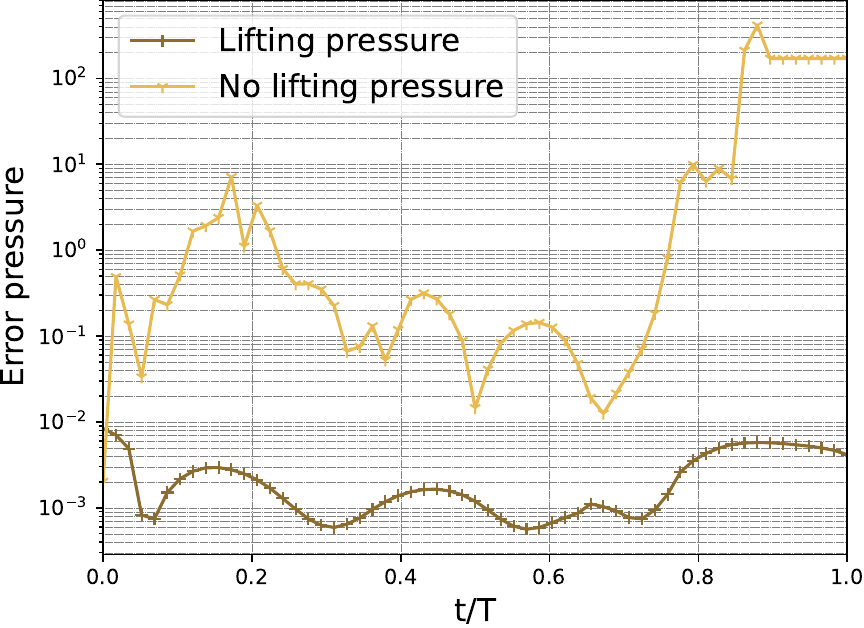}}\\
	\caption{Comparison of FOM-ROM errors for velocity and pressure with $N = 6$, both with and without the use of the lifting function for the pressure.}
	\label{errup_aorta_noliftP}
\end{figure}

Figure \ref{errup_aorta_noliftP} highlights the significance of incorporating the lifting function for pressure within our ROM framework. As shown in Figure \ref{fig:error_dtp_noliftP}, neglecting the lifting function for pressure results in an error exceeding $10^2$. Due to the coupling between velocity and pressure, this omission also leads to a twofold increase in velocity error, as observed in Figure \ref{fig:erru_aorta_noliftP}. The essential role of the lifting function for velocity has already been well established in the literature \cite{star2019novel}.

Table \ref{tab:nn} presents the hyperparameters of the neural network used to interpolate the outflow pressure at new time instances. Hyperparameter tuning is carried out by gradually increasing the learning rate and the number of hidden neurons for a fixed activation function and number of hidden layers, until a satisfactory decay of the loss function is observed, ensuring that overfitting is avoided. The dataset \eqref{time_set} is normalized and split into train ($80\%$ of the dataset) and test ($20\%$ of the dataset) sets.
\begin{table}
\centering
\renewcommand{\arraystretch}{1.5}
\caption{Hyperparameters used for the feedforward neural network.}
\begin{tabular}{ccccc}
\hline
\rowcolor{gray!20} Neurons per layer & Activation funtion  & Number of epochs & Learning rate & Hidden layers  \\
\hline
150  &   Softplus   & 50000 & 5$\cdot 10^{-6}$ & 2
%\multirow{2}*{147048 $\times$ 20} 
\\
\hline
\end{tabular}
\label{tab:nn}
\end{table}

\begin{figure}
	\centering
    \subfloat[][\label{fig:erru_aorta_newt}]{\includegraphics[width=.4\textwidth]{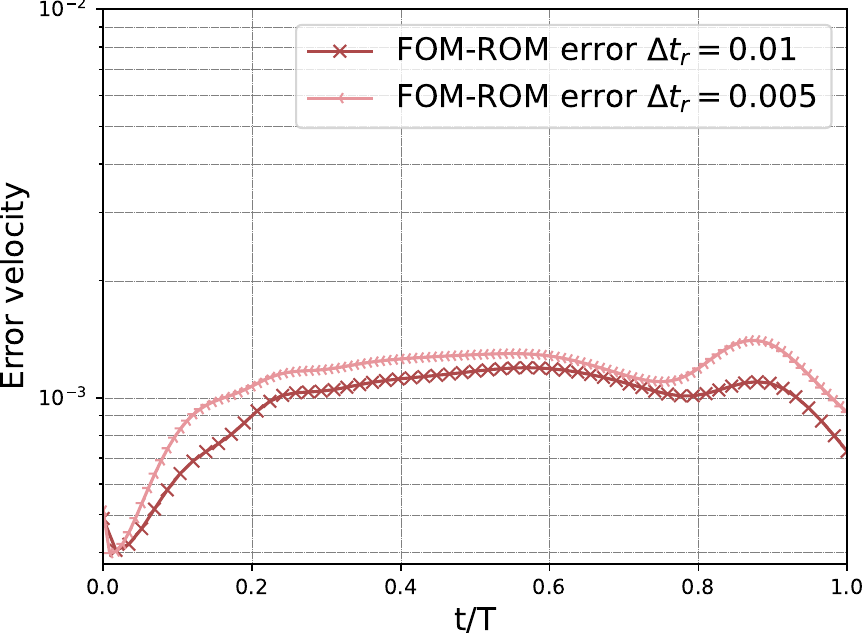}}\hspace{2ex}
 	\subfloat[][\label{fig:error_dtp_newt}]{\includegraphics[width=.4\textwidth]{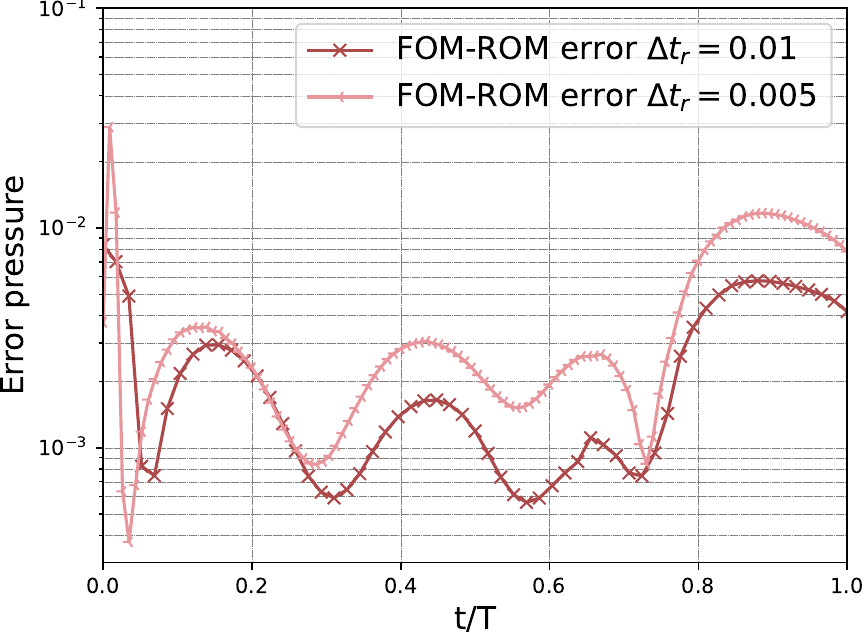}}\\
	\caption{Comparison of FOM-ROM errors for velocity and pressure with $N = 6$ for the two numerical experiments performed.}
    \label{errup_aorta_newt}
\end{figure}

Figure~\ref{errup_aorta_newt} presents a comparison of the FOM-ROM error for velocity and pressure when computing the reduced solution with $\Delta t_r = 0.01$ (as in the previous numerical experiment) and a refined time step of $\Delta t_r = 0.005$. As shown in Figure~\ref{fig:erru_aorta_newt}, reducing $\Delta t_r$ results in a slight increase in the velocity and pressure errors. However, the overall error trends remain consistent, demonstrating the robustness of our ROM framework. This confirms its reliability in predicting solution dynamics even at time instances where no full-order outflow pressure data has been stored.

\begin{figure}[!htb]
    \centering
    \vspace{3ex}
        \includegraphics[scale=0.28, trim = 0cm 0.85cm 0cm 0cm]{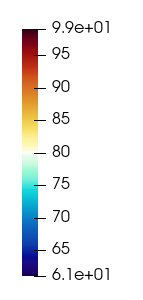}
        \hspace{-4ex}
        \begin{overpic}[width=0.25\textwidth]{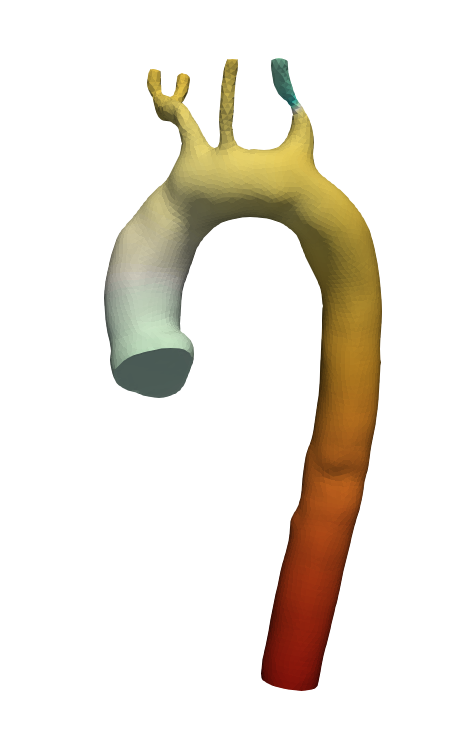} 
        \put(16,100){$p$ FOM}
        \put(-3,16){\rotatebox{90}{(mmHg)}}
        \end{overpic}
        \hspace{-5ex}
        \begin{overpic}[width=0.25\textwidth]{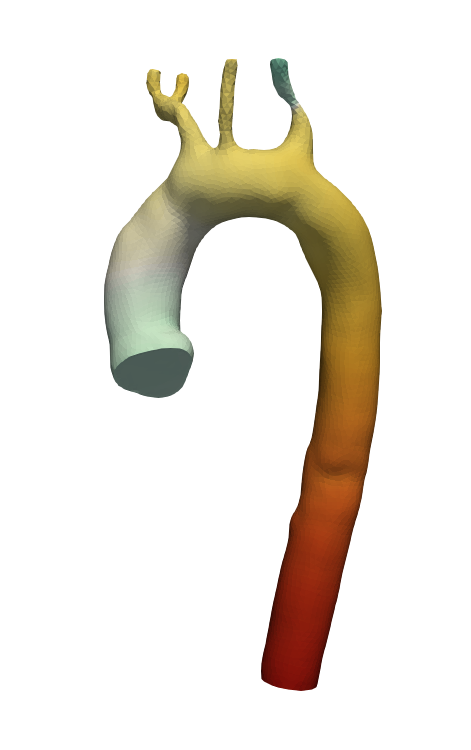} 
        \put(16,100){$p$ ROM}
        %\put(-10,107){$t=5.475$ s}%$t=5.5$ s}
        \end{overpic}
        \hspace{-1ex}
        \vline
        \hspace{-1ex}
        % \hspace{2ex}
        % \includegraphics[scale=0.28, trim = 0cm 0.85cm 0cm 0cm]{img/legend_p_aorta_1.png}
        \begin{overpic}[width=0.25\textwidth]{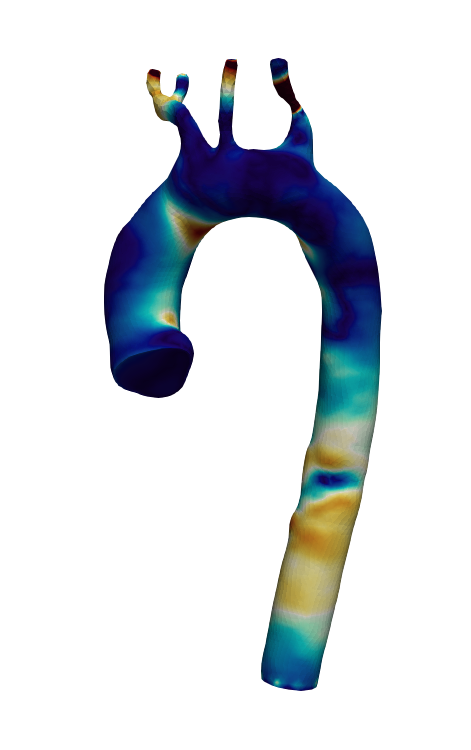} 
        \put(5,100){$\mid$ $p$ FOM - $p$ ROM $\mid$}
        %\put(65,16){\makebox(0,0)[r]{\rotatebox{90}{(mmHg)}}}
        %\put(-10,107){$t=5.5$ s}
        \end{overpic}
        %\hspace{2ex}
        \includegraphics[scale=0.28, trim = 0cm 0.85cm 0cm 0cm]{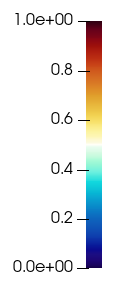}
            
\caption{Qualitative comparison of FOM-ROM pressure for the new time $t=5.475$ s.}
\label{fig:p-fom-rom-aorta-newt}
\end{figure}

Qualitative comparisons reveal a strong agreement in the pressure distribution between the FOM and the ROM, as illustrated in Figure~\ref{fig:p-fom-rom-aorta-newt} for the test value $t = 5.475$ s.

% \begin{figure}[!htb]
%     \centering
%         \raisebox{0.6cm}{\includegraphics[scale=0.25, trim = 0cm 0cm 0cm 0cm]{img/legend_p.png}}
%         \begin{overpic}[width=0.3\textwidth]{img/p_FOM_5475.png} 
%         \put(19,96){p FOM}
%         \put(-8,25){\rotatebox{90}{(mmHg)}}
%         \end{overpic}
%         \begin{overpic}[width=0.3\textwidth]{img/p_ROM_15.png} 
%         \put(19,96){p ROM}
%         %\put(-10,107){$t=5.475$ s}%$t=5.5$ s}
%         \end{overpic}
% \caption{Qualitative comparison of FOM-ROM pressure for the new time $t=5.475$ s.}
% \label{fig:p-fom-rom-aorta-newt}
% \end{figure}

%%%%%%%%%%%%%%%%%%%%%%%%%%%%%%%%%%%%%%%%%%%%%%%%%%
%%%%%%%%%%%%%%%%%%%%%%%%%%%%%%%%%%%%%%%%%%%%%%%%%%

\begin{figure}[!htb]
    \centering
    \includegraphics[scale=0.2, trim = 0cm 0.85cm 0cm 0cm]{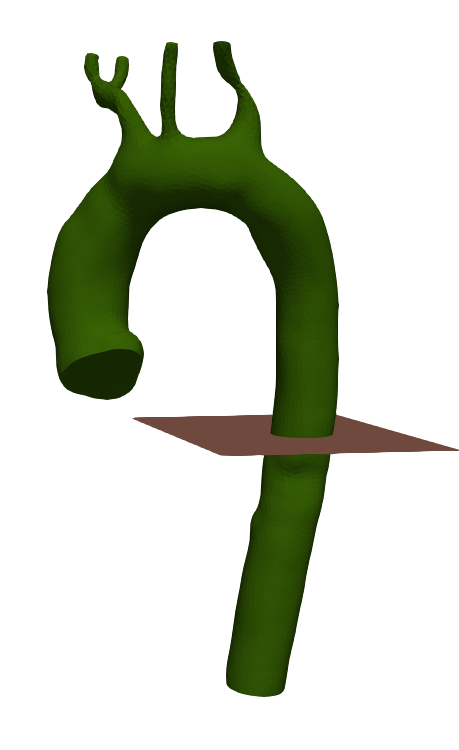}\\ 
        \raisebox{1.8cm}{\includegraphics[scale=0.25, trim = 0cm 0.85cm 0cm 0cm]{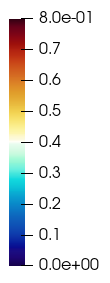}}
        \hspace{2ex}
        \begin{overpic}[width=0.3\textwidth]{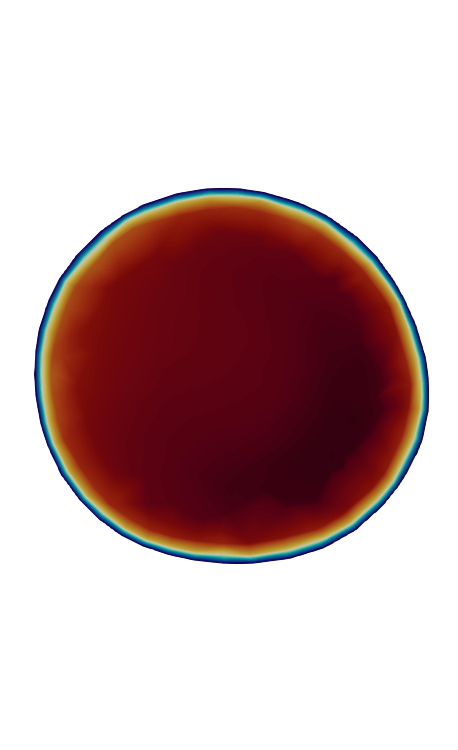} 
        \put(19,82){$\bm u$ FOM}
        \put(-12,45){\rotatebox{90}{(m/s)}}
        \end{overpic}
        \begin{overpic}[width=0.3\textwidth]{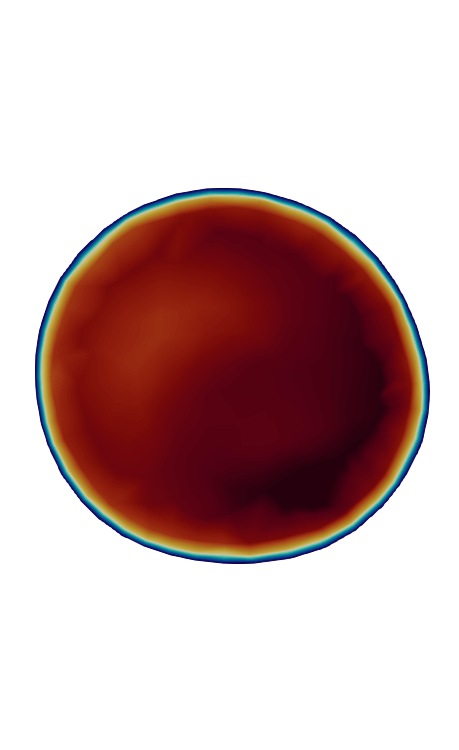} 
        \put(19,82){$\bm u$ ROM}
        %\put(-10,107){$t=5.475$ s}%$t=5.5$ s}
        \end{overpic}
\caption{Qualitative comparison of FOM-ROM velocity for the new time $t=5.955$ s on the slice (brown plane) shown in the descending aorta.}
\label{fig:u-fom-rom-aorta-newt}
\end{figure}

To further analyze the predicted solutions, Figure~\ref{fig:u-fom-rom-aorta-newt} presents the velocity magnitude of the FOM and ROM on a slice of the descending aorta. The ROM accurately captures the flow structures, exhibiting only minor discrepancies compared to the FOM.

To assess the computational efficiency, we analyze performance on an Intel(R) Core(TM) i7-7700 CPU @ 3.60GHz with 16GB RAM. A full-order model (FOM) simulation takes approximately one day, while the POD algorithm completes in 27 seconds, and neural network training requires 87 seconds. The evaluation phase is highly efficient, taking less than a second (approximately 285 milliseconds). Consequently, the achieved speedup -- defined as the ratio between the computational time of the FOM simulation and that of solving the reduced-order system -- is substantial, reaching the order of $\mathcal{O}(10^{5})$.

\section{Conclusions}

\textcolor{black}{In this work, we have presented a hybrid ROM approach for cardiovascular flows, where time is the only varying parameter. Non-homogeneous boundary conditions, which are not automatically preserved at the ROM level, are treated with the lifting function method \cite{stabile2017pod,graham1999optimal,gunzburger2007reduced}, which shifts high-fidelity solutions. Boundary conditions are then introduced during the reconstruction phase, allowing the ROM to incorporate arbitrary boundary values.} 

\textcolor{black}{While this technique is well established for the velocity field in CFD applications, its extension to pressure, with realistic non-homogeneous outlet conditions, remains unexplored. Here, the Windkessel model is used to provide outflow pressures at the FOM level, and neural networks are employed within the ROM to approximate these values efficiently, enabling online evaluation beyond the original FOM snapshots.}

\textcolor{black}{The ROM framework has been tested on a 3D patient-specific aortic arch, and the ROM demonstrates good accuracy in terms of both quantitative errors and qualitative comparison of pressure and velocity fields. The computational speed-up is significant, requiring less than one second for ROM evaluations once the FOM is trained.} %These results suggest that the methodology provides a solid foundation for multiparametric cardiovascular models and shows promising generalization capabilities even with a limited training dataset.  

\textcolor{black}{Extending the approach to multiparametric settings, including physical and geometric variations like in \cite{Siena2022,siena2023fast,Stabile2020,Padula2024}, represents a natural next step. Machine learning techniques, such as autoencoders \cite{phillips2021autoencoder}, could further enhance the ROM by capturing nonlinear features in high-fidelity solutions more efficiently.}

\textcolor{black}{An interesting future direction concerns the parametrization and calibration of the Windkessel model using clinical data, leveraging the ROM framework to efficiently %identify the most representative parameter configurations and assess the physiological validity of 
assess the model in patient-specific settings. The use of fully nonlinear data-driven ROMs is also a very promising avenue; for instance, approaches based on autoencoders could be employed to capture nonlinear dependencies more effectively.}

\bibliographystyle{plainnat}
\bibliography{bib.bib}

\end{document}